\documentclass[12pt,oneside]{amsart}
\usepackage{amssymb,amscd}
\usepackage{url}
\usepackage{mathtools}
\usepackage{enumerate}

\usepackage{cite, verbatim}

\newtheorem{theorem}{Theorem}[section]
\newtheorem{lemma}[theorem]{Lemma}
\newtheorem{proposition}[theorem]{Proposition}
\newtheorem{corollary}[theorem]{Corollary}

\theoremstyle{definition}

\theoremstyle{remark}
\newtheorem{remark}[theorem]{Remark}

\DeclareMathOperator{\Sym}{Sym}
\DeclareMathOperator{\Alt}{Alt}
\DeclareMathOperator{\Aut}{Aut}
\DeclareMathOperator{\AGL}{A\Gamma L}
\DeclareMathOperator{\GamL}{\Gamma L}
\DeclareMathOperator{\GamO}{\Gamma O}
\DeclareMathOperator{\GL}{GL}
\DeclareMathOperator{\SL}{SL}

\DeclareMathOperator{\Sz}{Sz}
\DeclareMathOperator{\GF}{GF}
\DeclareMathOperator{\VO}{VO}
\DeclareMathOperator{\VSz}{VSz}
\DeclareMathOperator{\VD}{VD}
\DeclareMathOperator{\ord}{ord}
\DeclareMathOperator{\Soc}{Soc}

\title[The automorphism groups of small affine rank 3 graphs]{The automorphism groups of small affine rank 3 graphs}

\author{Jin Guo}
\address{School of Mathematics and Statistics, Hainan University, Haikou 570228, Hainan, P. R. China}
\email{jguo@hainanu.edu.cn}

\author{Andrey V. Vasil'ev}
\address{Sobolev Institute of Mathematics, Novosibirsk 630090, Russia}
\email{vasand@math.nsc.ru}

\author{Rui Wang}
\address{School of Mathematics and Statistics, Hainan University, Haikou 570228, Hainan, P. R. China}
\email{wangr15155615917@163.com}

\begin{document}

\begin{abstract}
A \emph{rank $3$ graph} is an orbital graph of a rank $3$ permutation group of even order. Despite the classification of rank 3 graphs being complete, see, e.g., Chapter 11 of the recent monograph \emph{Strongly regular graphs} by Brouwer and Van~Maldeghem, the full automorphism groups of these graphs (equivalently, the $2$-closures of rank $3$ groups) have not been explicitly described, though a lot of information on this subject is available. In the present note, we address this problem for the affine rank 3 graphs. We find the automorphism groups for finitely many relatively small graphs and show that modulo known results, this provides the full description of the automorphism groups of the affine rank 3 graphs, thus reducing the general problem to the case when the socle of the automorphism group is nonabelian simple.

\smallskip

\noindent\textsc{Keywords:} affine permutation group, closure of permutation group, orbital graph, full automorphism group of graph, rank 3 group, rank 3 graph.
\smallskip

\noindent\textsc{MSC:} 20B25, 05E18, 05E30.
\end{abstract}

\maketitle

\section{Introduction}\label{s:intro}

Let $G\leq\Sym(\Omega)$ be a permutation group on a finite set $\Omega$. The orbits in the induced action of $G$
on~$\Omega \times \Omega$ are called 2-\textit{orbits} or {\em orbitals} of~$G$. The {\em $2$-closure} $G^{(2)}$ of $G$ is the largest permutation group having the same 2-orbits as~$G$. A graph $\Gamma$ with the vertex set $\Omega$ and the arc set $E$ equal to one of the orbitals of $G$ is called an {\em orbital graph} of $G$.

The {\em rank} of $G$ is the number of $2$-orbits of $G$. If the rank of $G$ equals $3$, then $G$ is called a {\em rank $3$ group}. Let $D,E,F$ be the orbitals of a rank $3$ group $G$. Clearly, $G$ must be transitive, so one of its orbitals, say $D$, is the diagonal of $\Omega\times\Omega$. If the order of $G$ is odd, then two other orbitals $E$ and $F$ are transpose to each other, so the corresponding orbital graphs are {\em tournaments}, i.e., complete directed graphs, opposite to each other. If the order of $G$ is even, then $E$ and $F$ are symmetric, and the corresponding two graphs are complementary undirected graphs called {\em rank $3$ graphs} (corresponding to a group $G$). As readily seen, $G$ is an automorphism group of such a graph~$\Gamma$ and $G^{(2)}=\Aut(\Gamma)$ is the (full) automorphism group of it (the automorphism groups are obviously the same for complementary graphs).

It is clear that a rank $3$ graph is {\em strongly regular}, that is the number of the common neighbours of two distinct vertices of the graph depends only on whether they are equal, adjacent, or nonadjacent. The rank $3$ graphs form one the most important and well-studied subclass of the class of strongly regular graphs. Nevertheless, there are still challenging problems concerning these graphs. For example, no explicit combinatorial characterization of the rank $3$ graphs is known, so one cannot effectively (i.e., by a polynomial-time algorithm) check whether a given strongly regular graph is a rank $3$ graph or not.

The classification of the rank $3$ graphs and groups is due to Banai, Higman, Foulser, Kallaher, Kantor, Liebler, Liebeck, Saxl and many others. Chapter 11 of the recent monograph \emph{Strongly regular graphs} by Brouwer and Van~Maldeghem \cite{srgw} contains the description of all pairs $(\Gamma, G)$ with $\Gamma$ a rank $3$ graph and $G$ an automorphism group of this graph acting rank~$3$. It is surprising that despite this, the full automorphism groups of rank $3$ graphs (equivalently, the $2$-closures of rank $3$ groups) are apparently not completely described. Undoubtedly, a lot of information on this matter is known, see again \cite[Chapter~11]{srgw}. In particular, for every rank $3$ graph $\Gamma$ of degree at most $1024$, the automorphism group of $\Gamma$ is given explicitly (see \cite[Table~11.8]{srgw}).

Significant progress in the description of the $2$-closures of rank $3$ groups was recently obtained by Skresanov in~\cite[Theorems~1.1 and 1.2]{SkrARS}, see also Theorems~\ref{class} and~\ref{a1inv} in Section~\ref{s:known} below. These results allow him to construct a polynomial-time algorithm for finding the $2$-closure of a rank $3$ permutation group \cite{SkrJA}. Equivalently, this algorithm finds the automorphism group of a rank $3$ graph (if some group of rank $3$, which this graph corresponds to, is given as input to the algorithm). The ultimate goal of \cite{SkrARS, SkrJA} --- to obtain a computational algorithm --- left some classes of small affine groups of rank 3 without close consideration. In the present paper, we fill this gap --- see Propositions~\ref{p:onedim}--\ref{p:exceptC} below --- by giving the description of 2-closures of the affine rank~3 groups in question. This reduces the description of the automorphism groups of rank $3$ graphs to the case when the socle of the automorphism group is nonabelian simple.

Recall that a permutation group $G \leq \Sym(\Omega)$ is called {\em affine}, if it has a regular normal elementary abelian $p$-subgroup $V$. Then $G =V:G_0$ is a split extension of $V$ and a point stabilizer $G_0$. The set $\Omega$ can be identified with $V$ considered as a vector space over the prime field $\GF(p)$ in such a way that $V$ itself acts on this space by translations, and $G_0$ being the stabilizer of the zero vector acts on it as a subgroup of $\GL(V)$. Such an identification provides a natural embedding $G\leq\operatorname{AGL}(V)$.

Together with the known results (see Section~\ref{s:known}), our efforts in Section~\ref{s:proof} allow us to claim that the following general assertion holds.

\begin{theorem}\label{t:2clos}
If $G$ is a finite affine rank $3$ group, then the $2$-closure $G^{(2)}$ is known.
\end{theorem}

A graph $\Gamma$ is called an {\em affine rank $3$ graph} if its vertex set can be identified with a vector space $V=GF(p)^m$ and $\Aut(\Gamma)\cap\operatorname{AGL}(V)$ is an affine rank~$3$ group. In other words, the edge set of $\Gamma$ coincides with one of irreflexive orbitals of a rank $3$ group $G\leq\operatorname{AGL}(V)$. In this notation, Theorem~\ref{t:2clos} implies that for every affine rank $3$ graph $\Gamma$ given by some affine rank $3$ group, the automorphism group $\Aut(\Gamma)$ is known.  Modulo the known results (see the last paragraph of Section~\ref{s:known}), it follows that the next statement holds.

\begin{corollary}\label{c:aut}
Suppose that $\Gamma$ is a finite rank $3$ graph. Then either the automorphism group $\Aut(\Gamma)$ is known or the socle of $\Aut(\Gamma)$ is a known nonabelian simple group.
\end{corollary}

Our statement of Corollary~\ref{c:aut} leaves open the question on the full automorphism group of a rank $3$ graph, having a simple socle. In fact, the classification of the rank $3$ groups in~\cite{srgw} with a little help of the description of the maximal subgroups of finite groups of Lie type in~\cite{bray,KL,Craven} provide all necessary information to address this question. Nevertheless, from our point of view, it would be useful to have the explicit description of the full automorphism groups of the rank 3 graphs in this last case.

\section{The known facts and plan of the proof}\label{s:known}

Our notation is mostly standard.

For a prime $p$ and natural number $a$, the multiplicative order of $a$ in the prime field $GF(p)$ is denoted by $\ord_p(a)$.

Given groups $A$ and $B$, $A\times B$, $A\circ B$, and $A:B$ stand for their direct, central and semidirect products, respectively. The $m$th direct power of $A$ is~$A^m$. The cyclic group of order $n$ is denoted simply by $n$ except the cases where this might confuse a reader. Along this way, $p^m$ stands for the elementary abelian group of order $p^m$. The groups $D_{n}$ and $Q_8$ are the dihedral group of order $n$ and the quaternion group, respectively. Simple sporadic groups and simple groups of Lie type are denoted according to \cite{Atlas}, while the notation of other finite groups of Lie type follows \cite{KL}. We refer to the symmetric and alternating groups of degree~$n$ as $\Sym(n)$ and $\Alt(n)$, respectively.

If a permutation group $G\leq\Sym(\Omega)$ is transitive, the orbitals of $G$ are in the natural one-to-one correspondence with the orbits of a point stabilizer $G_\alpha$, in particular, the diagonal orbital maps to the singleton orbit $\{\alpha\}$. The sizes of the other orbits (corresponding to irreflexive orbitals) are called the {\em subdegrees} of~$G$. Thus, a rank 3 group has exactly two subdegrees.

If $A\leq\Sym(\Gamma)$ and $B\leq\Sym(\Delta)$, then $A\wr B\leq\Sym(\Gamma\times\Delta)$ and $A\uparrow B\leq\Sym(\Gamma^\Delta)$ denote their imprimitive and primitive wreath products, respectively.

The affine group $G=V:G_0$ is primitive if and only if $G_0\leq\GL(V)$ is irreducible. In this case, $V=\Soc(G)$ is a unique minimal normal subgroup of~$G$.

If the zero stabilizer $G_0$ acts semilinearly on $V$ as an $F$-vector space, where $F=\GF(q)$ and $q$ is a power of~$p$,
then we write $G_0 \leq \GamL_m(q)$, where $\GamL_m(q)$ is the full semilinear group, and identify $V$ with~$F^m$. We also refer to $\AGL_m(q)=V:\GamL_m(q)$ as the full affine semilinear group. If the field is clear from the context, we may use $\GamL(V) = \GamL_m(q)$ and $\AGL(V) = \AGL_m(q)$ instead.

An affine rank 3 graph $\Gamma$ is {\em$m$-dimensional}, if the point set of $\Gamma$ can be identified with a vector space $V=F^m$ in such a way that $\Aut(\Gamma)\cap\AGL(V)$ is a rank $3$ group. In this case, the edge set of $\Gamma$ coincides with one of irreflexive orbitals of an affine semilinear rank $3$ group $G\leq\AGL_m(q)$. Other notations for rank 3 graphs are taken from \cite{srgw}.

We start with the main result of \cite{SkrARS}.

\begin{theorem}{\rm\cite[Theorem 1.1]{SkrARS}}\label{class}
	Let $G$ be a rank~$3$ permutation group on a finite set $\Omega$. Then either $G$ is one of the groups from Table~{\rm\ref{exceptions}}, or exactly one of the following is true.
	\begin{enumerate}[{\rm(i)}]
		\item $G$ is imprimitive, i.e.\ it preserves a nontrivial decomposition
			$\Omega = \Delta \times X$.
			Then $$G^{(2)} = \Sym(\Delta) \wr \Sym(X).$$
		\item $G$ is primitive and preserves a product decomposition
			$\Omega = \Delta^2$. Then $$G^{(2)} = \Sym(\Delta) \uparrow \Sym(2).$$
		\item $ G $ is primitive almost simple with socle $ L $, i.e.\ $ L \unlhd G \leq \Aut(L) $.
			Then $ G^{(2)} = N_{\Sym(\Omega)}(L) $, and $ G^{(2)} $ is almost simple with socle $ L $.
		\item $ G $ is a primitive affine group, i.e.\ $ G \leq \AGL_a(q) $ for some $ a \geq 1 $
			and a prime power~$ q $, moreover, $ G $ does not stabilize a product decomposition.
			Set $ F = \GF(q) $. Then $ G^{(2)} $ is also an affine group and exactly one of the following holds.
			\begin{enumerate}[{\rm(a)}]
				\item $ G \leq \AGL_1(q) $. Then $ G^{(2)} \leq \AGL_1(q) $.
				\item $ G \leq \AGL_{2m}(q) $ preserves the bilinear forms graph $ H_q(2, m) $, $ m \geq 3 $. Then
					$$G^{(2)} = F^{2m} : ((\GL_2(q) \circ \GL_m(q)) : \Aut(F)).$$
				\item $ G \leq \AGL_{2m}(q) $ preserves the affine polar graph $ \VO_{2m}^{\epsilon}(q) $,
					$ m \geq 2 $, $ \epsilon = \pm $. Then
					$$G^{(2)} = F^{2m} : \GamO^{\epsilon}_{2m}(q).$$
				\item $ G \leq \AGL_{10}(q) $ preserves the alternating forms graph $ A(5, q) $. Then
					$$G^{(2)} = F^{10}: ((\GamL_5(q) / \{ \pm 1 \}) \times (F^* / (F^*)^2)).$$
				\item $ G \leq \AGL_{16}(q) $ preserves the affine half spin graph $ \VD_{5,5}(q) $.
					Then $ G^{(2)} \leq \AGL_{16}(q) $ and we have
					$$G^{(2)} = F^{16}:((F^* \circ \operatorname{Inndiag}(D_5(q))):\Aut(F)).$$
				\item $ G \leq \AGL_4(q) $ preserves the Suzuki-Tits ovoid graph $ \VSz(q) $,
					$ q = 2^{2e+1} $, $ e \geq 1 $. Then
					$$ G^{(2)} = F^4:((F^* \times \Sz(q)):\Aut(F)).$$
			\end{enumerate}
	\end{enumerate}
\end{theorem}

\begin{table}[ht]
\caption{Possible exceptions to Theorem~\ref{class}}\label{exceptions}
\begin{tabular}{l l l}
\hline
type of a group &  appearance & parameters\\
\hline

	 Almost simple & Corollary~\ref{c:aut} & Table~\ref{exceptAS}\\
	$G\leq\AGL_1(q)$ & Prop.~\ref{p:onedim} &  Table~\ref{a1exctab} and Table~\ref{not1dim}\\
	Extraspecial class (B) & Prop.~\ref{p:exceptB} & Table~\ref{exceptB} \\
    Exceptional class (C) & Prop.~\ref{p:exceptC} & Table~\ref{exceptC} \\
		\hline
\end{tabular}
\end{table}

\begin{remark} According to  Table~\ref{exceptions}, all possible exceptions to Theorem~\ref{class} are divided in four types (see the first column of the table). The second column contains the references to the assertions of the present paper, where the groups of the corresponding type are considered. The third column refers to the tables, where the parameters of the groups are listed. For $G\leq\AGL_1(q)$, we give the reference to two tables: Table~\ref{a1exctab} is~\cite[Table~7]{SkrARS} --- we use it in the proof of Proposition~\ref{p:onedim}, while Table~\ref{not1dim} contains the results of our considerations.
\end{remark}

The classification of the affine rank 3 groups was completed by Liebeck in \cite{liebeckAffine} (see also \cite[Theorem~11.4.1]{srgw}). He divided the affine groups of rank 3 into three classes: (A), (B), and~(C). Class (A) includes all the infinite series, while `extraspecial' class (B) and `exceptional' class (C) consist only of finitely many groups.  It follows from Theorem~\ref{class} and Table~\ref{exceptions} that for an affine rank 3 group $G$, the precise structure of $G^{(2)}$ is  unknown only if $G\leq\AGL_1(q)$ for class~(A), or $G$ belongs to one of the classes (B)  and~(C). The subclass (A1) of the one-dimensional affine groups $G\leq\AGL_1(q)$ was also considered in \cite{SkrARS} for sufficiently large~$q$.

\begin{theorem}{\rm\cite[Theorem 1.2]{SkrARS}}\label{a1inv}
	Let $G$ be a primitive affine permutation group of rank~$3$
	and suppose that $G \leq \AGL_1(q)$ for some prime power $q$.
	Then either $G^{(2)}$ lies in $\AGL_1(q)$, or degree and
	the smallest subdegree of $G$ are listed in Table~{\rm\ref{a1exctab}}.
\end{theorem}

\begin{table}[ht]
\caption{Possible exceptions to Theorem~\ref{a1inv}}\label{a1exctab}
	(A)
\begin{tabular}{c|c|c|c|c|c|c|c|c|c|c}
	\hline
	degree & $2^4$ & $ 2^6 $ & $ 2^8 $ & $ 2^{10} $ & $ 2^{12} $ & $ 2^{16} $ & $3^2$ & $3^4$ & $3^6$ & $ 5^2 $\\
	subdegree & 5 & 21 & 51 & 93 & 315 & 3855 & 4 & 16 & 104 & 8\\
	\hline
\end{tabular}

	(B)
\begin{tabular}{c|c|c|c|c|c|c|c}
	\hline
	degree & $ 3^4 $ & $ 3^6 $ & $ 7^2 $ & $ 7^4 $ & $ 17^2 $ & $ 23^2 $ & $ 47^2 $\\
	subdegree & 16 & 104 & 24 & 480 & 96 & 264 & 1104\\
	\hline
\end{tabular}

	(C)
\begin{tabular}{c|c|c}
	\hline
	degree & $ 3^4 $ & $ 89^2 $\\
	subdegree & 40 & 2640\\
	\hline
\end{tabular}
\end{table}

\begin{remark} Table~\ref{a1exctab} repeats \cite[Table~7]{SkrARS}. It contains three rows noted by (A), (B), and~(C). According to \cite[Lemmas~3.3-3.5]{SkrARS}, a rank 3 group $G$ that does not embed into $\AGL_1(q)$ and has the same parameters: degree and subdegree, as in the table, might belong only to class (A), (B), or~(C) depending on which row of the table the parameters were taken.
\end{remark}

The first step in order to prove Theorem~\ref{t:2clos} is to find the $2$-closures of subgroups of $\AGL_1(q)$ with the parameters listed in Table~\ref{a1exctab}. Observe that the most of these groups have degree less than or equal to $1024$, so their $2$-closures as the automorphism groups of the corresponding rank 3 graphs can be extracted from~\cite[Table~11.8]{srgw}. These groups alongside with the larger one-dimensional affine groups with parameters from Table~\ref{a1exctab} are considered in Proposition~\ref{p:onedim}. Resulting Table~\ref{not1dim} collects the information on all groups $G\leq\AGL_1(q)$ with $G^{(2)}\not\leq\AGL_1(q)$.

Below we provide the information on the one-dimensional affine rank 3 groups and graphs that we need in the proof. These groups was described by Foulser and Kallaher \cite{foulserRank3}. Here we follow to more recent papers by Peisert \cite{peisert}  and especially by Muzychuk \cite{muzychukOneDim}. In the latter paper, the one-dimensional affine rank 3 graphs were classified in the following sense. If $G\leq\AGL_1(q)$ is a one-dimensional affine rank 3 group, then $G^{(2)}\cap\AGL_1(q)$ is described.

Let $p$ be a prime, $q=p^d$, $F=GF(q)$, $\omega\in F$ a generator of the multiplicative group $F^*$, and $\Phi=\Aut(F)=\langle\phi\rangle$, where $\phi: F\rightarrow F$ is a Frobenius automorphism of $F$ that moves each $x$ in $F$ to $x^p$. Then $\AGL_1(q)\simeq F:(F^*:\Phi)$ and $G\leq \AGL_1(q)$ is isomorphic to $F:G_0$, where the zero stabilizer $G_0$ is a subgroup in $\GamL_1(q)$. Denote by $\widehat\omega$ and $\widehat\phi$ the elements in $\GamL_1(q)$ corresponding to $\omega$ and $\phi$ under the natural isomorphism $\GamL_1(q)\simeq F^*:\Phi$.

\begin{lemma}{\rm\cite[Propositon~1]{muzychukOneDim}}\label{VLS}
Let $p$ and $e$ be two primes satisfying $\ord_e(p)=e-1$ and $q=p^{k(e-1)}$, where $k$ is a positive integer. Let $C=\langle\omega^e\rangle\leq F^*$. Then the subgroup $G_0=\langle\widehat\omega^e,\widehat\phi\rangle\leq\GamL_1(q)$ of order $k(e-1)(q-1)/e$ has a cyclic subgroup $\widehat{C}=G_0\cap\GL_1(q)=\langle\widehat\omega^e\rangle$ of order $(q-1)/e$. The group $G_0$ has two orbits on $F^*$ of sizes $(q-1)/e$ and $(q-1)(e-1)/e$, namely, $C$ and $F\setminus C$. If $H\leq\GamL_1(q)$ has the same orbits, then $(F:H)^{(2)}\cap\AGL_1(q)=F:G_0$.
\end{lemma}

The corresponding rank 3 graph $\Gamma$ is the {\em Van Lint--Schrijver graph} if $e>2$ and the {\em Paley graph} if $e=2$. If $p\equiv3\pmod4$ and $k$ is even, then there is another rank 3 graph $\Gamma$ on $q=p^k$ vertices with the same orbit lengths as the Paley graph. This graph is called the {Peisert graph} and is defined below.

\begin{lemma}{\rm\cite[Propositon~2]{muzychukOneDim}}\label{Peisert}
Let $p\equiv3\pmod4$ be a prime and $k$ an even positive integer. Let $C=\langle\omega^4\rangle$. Then the subgroup $G_0=\langle\widehat\omega^4,\widehat\phi\widehat\omega\rangle\leq\GamL_1(q)$ of order $k(q-1)/4$ has a cyclic subgroup $\widehat{C}=G_0\cap\GL_1(q)=\langle\widehat\omega^e\rangle$ of order $(q-1)/4$. The group $G_0$ has two orbits on $F^*$ of the same size $(q-1)/2$, namely, $C\cup C\omega$ and $C\omega^2\cup C\omega^3$. If $H\leq\GamL_1(q)$ has the same orbits, then $(F:H)^{(2)}\cap\AGL_1(q)=F:G_0$.
\end{lemma}

Observe that in the notation of Lemma~\ref{Peisert}, the two orbits $C\cup C\omega^2$ and $C\omega\cup C\omega^3$ of the subgroup $\langle\widehat\omega^2,\widehat\phi\rangle\leq\GamL_1(q)$ yields the complementary Paley graphs. It is also worth mentioning that for $q=3^2$, there is an isomorphism between the corresponding Paley and Peisert graphs, while for $q>3^2$, they are always nonisomorphic, see \cite{peisert}. Describing the finite self-complementary symmetric graphs, Peisert proved \cite{peisert} that besides the one-dimensional Paley and Peisert graphs (cf. Lemmas~\ref{VLS} and~\ref{Peisert}) there are three self-complementary symmetric graphs on $49$, $81$ and $529$ vertices, which are orbital graphs of 2-dimensional affine groups. Two of these graphs, on $49$ and $81$ vertices, are isomorphic to the one-dimensional Peisert graphs, though their full automorphism groups are 2-dimensional affine groups. The third graph on $529=23^2$ vertices, it is called the {\em sporadic Peisert graph}, cannot be constructed as an orbital graph of one-dimensional affine group. In fact, the full automorphism group of this graph belongs to the `extraspecial' class (B), see below. We summarize the information from this paragraph as follows.

\begin{lemma}{\rm\cite[Theorem 7.1]{peisert}}\label{PP}
The Paley and Peisert graphs together with the sporadic Peisert graph on $23^2$ vertices are the only self-complementary symmetric graphs. They are not isomorphic to each other except when degree $q=3^2$. The automorphism groups of these graphs are known, they are subgroups of $\AGL_1(q)$ except when $q\in\{7^2,3^4,23^2\}$.
\end{lemma}

It turns out that each one-dimensional affine rank 3 graph is an orbital graph of one of the rank 3 groups described above.

\begin{theorem}{\rm\cite[Theorem~3]{muzychukOneDim}}\label{Muz}
Let $\Gamma$ be a one-dimensional affine rank $3$ graph. Then, up to a complement, $\Gamma$ is either the
Van Lint--Schrijver, the Paley or the Peisert graph.
\end{theorem}
\medskip

As mentioned above, the rest of the affine groups of rank 3 that we need to consider are from the classes (B) and (C) (see Tables~\ref{exceptB} and~\ref{exceptC}). For most of them  their $2$-closures can be extracted from~\cite[Chapters~10,11]{srgw}. We deal with the groups from classes (B) and (C) in Propositions~\ref{p:exceptB} and~\ref{p:exceptC} and list the information on their $2$-closures in Tables~\ref{exceptB} and~\ref{exceptC}, respectively, thus completing the proof of Theorem~\ref{t:2clos}. Below we provide the necessary facts on these groups, see  \cite{srgw,liebeckAffine} for details.

Let $V=GF(p)^d$ be the socle of an affine primitive rank 3 group $G$ of degree $n=p^d$. The zero stabilizer $G_0$ is a subgroup in $\GL(V)\simeq\GL_d(p)$.

Extraspecial class (B). Here the zero stabilizer $G_0$ is contained in the normalizer $N=N_{GL(V)}(R)$ in $\GL(V)$ of some $r$-subgroup $R$ irreducible on $V$. In the most cases, $R$ is an extraspecial group of order $r^{1+2m}$, where $r\in\{2,3\}$ and $m\in\{1,2,3\}$. More precisely, $r=3$ only if $n=2^6$ and $m=1$. If $r=2$, then $R=R_i^m$, where $m\in\{1,2,3\}$, $R_1^m$ is the extraspecial group of type $2_{+}^{1+2m}$, $R_2^m$ is the extraspecial group of type $2_{-}^{1+2m}$, and $R_3^m$ is the central product of the cyclic group of order  $4$ and the extraspecial group of order $2^{1+2m}$. In fact, the normalizer $N$ may include more than one rank 3 subgroup, so there is some ambiguity in the description of these subgroups in Liebeck's classification. However, the maximal rank 3 subgroups from this class can be extracted from the proofs in \cite{foulserLowRank,liebeckAffine}. We list the information on `extraspecial' rank 3 groups in Table~\ref{exceptB} in the next section together with the description of their $2$-closures.

Exceptional class (C). Here the socle of $G_0/Z(G_0)$ is a simple group. We list the information on rank 3 groups from this class in Table~\ref{exceptC} in the next section together with the description of their $2$-closures. In most cases, the information on the 2-closures readily follows from \cite{srgw}, so in the proof of Proposition~\ref{p:exceptC}, we restrict ourselves to few necessary comments.
\medskip

We finish this section by showing how Corollary~\ref{c:aut} follows from Theorem~\ref{t:2clos} and known results.

\textit{Proof of Corollary~{\rm\ref{c:aut}}}.
Recall that a graph $\Gamma$ of rank $3$ is considered given if some group $G$ of rank $3$ such that $\Gamma$ is an orbital graph of~$G$ is given. Therefore, if $G$ is affine, then $\Aut(\Gamma)=G^{(2)}$ is known by virtue of Theorem~\ref{t:2clos}. If $G$ is imprimitive or preserves a product decomposition, then $\Aut(\Gamma)$ is also known, see Items (i) and (ii) of Theorem~\ref{class}. It follows that $G$ is an almost simple group, i.e., the socle $G$ is a nonabelian simple group. In the last case, all the situations, where the socles of $G$ and $G^{(2)}$ are distinct, are classified in~\cite[Theorem~1]{liebeck2Closure}. Applying also \cite[Propositions~1,2]{liebeck2Closure}, it is easy to find which of these groups are of rank~3. For convenience, we list this information in Table~\ref{exceptAS} below. Therefore, we come to Item~(iii) of Theorem~\ref{class}, thus completing the proof of Corollary~\ref{c:aut}.

\begin{table}[ht]
\caption{Almost simple groups of rank 3 with $\Soc(G)\neq\Soc(G^{(2)})$}\label{exceptAS}
\begin{tabular}{l|l|l|l}
\hline
$\Soc(G)$ & $\Soc(G^{(2)})$ & degree & conditions\\
\hline
	 $L_2(8)$ & $\Alt(9)$ & 36 & $G=\operatorname{P\Gamma{L}}_2(8)$\\
\hline
	 $M_{11}$ & $\Alt(11)$ & 55 & \\
\hline
	 $M_{12}$ & $\Alt(12)$ & 66 & \\
\hline
	 $M_{23}$ & $\Alt(23)$ & 253 & \\
\hline
	 $M_{24}$ & $\Alt(24)$ & 276 & \\
\hline
	 $\Alt(9)$ & $O^+_8(2)$ & 120\\
\hline
	 $G_2(q)$ & $O_7(q)$ & $q^3(q^3-1)/2$ & $q\in\{3,4,8\}$\\
\hline
	 $O_7(q)$ & $O_8^+(q)$ & $q^3(q^4-1)/\gcd(2,q-1)$ & $q\in\{2,3\}$\\
 		\hline
\end{tabular}
\end{table}

\section{Proof of Theorem~\ref{t:2clos}}\label{s:proof}

\begin{proposition}\label{p:onedim}
If  $G\leq\AGL_1(q)$, then either $G^{(2)}\leq\AGL_1(q)$ or $G^{(2)}$ is one of the groups in Table~{\rm\ref{not1dim}}.
\end{proposition}

\begin{table}[ht]
\caption{Groups $G\leq\AGL_1(q)$ with $G^{(2)}\not\leq\AGL_1(q)$}\label{not1dim}
\begin{tabular}{l|l|l|l|l}
\hline
$q=p^d$ & subdegs & $\Gamma$ & $G^{(2)}$ & refs\\
\hline

	$16=2^4$ & $5,10$ & VLS, $\VO^-_4(2)$ & $2^4:\Sym(5)$ & \cite[Table~11.8]{srgw}\\
\hline
	$25=5^2$ & $8,16$ & VLS, $5\times5$ & $(\Sym(5)\times\Sym(5)):2$ & \cite[Table~11.8]{srgw}\\
\hline
	$49=7^2$ & $24,24$ & Peisert, (B) & $7^2:(3\times\ SL_2(3))$ & \cite[\S~10.18]{srgw}\\
\hline
	$64=2^6$ & $21,42$ & VLS, $H_2(2,3)$ & $2^6:(\Sym(3)\times L_3(2))$ & \cite[Table~11.8]{srgw}\\
\hline
	$81=3^4$ & $16,64$ & VLS, $9\times9$ & $(\Sym(9)\times\Sym(9)):2$ & \cite[Table~11.8]{srgw}\\
\hline
	$81=3^4$ & $40,40$ & Peisert, (C) & $3^4:(\SL_2(5):2^2)$ & \cite[\S~10.30]{srgw}\\
\hline
	$256=2^8$ & $51,204$ & VLS, $\VO_4^-(4)$ & $2^8:(3\times \SL_2(2^4)):4$ & \cite[Table~11.8]{srgw}\\
\hline
	$729=3^6$ & $104,624$ & VLS, $H_3(2,3)$ & $3^6:(L_3(3)\times\GL_2(3))$ & \cite[Table~11.8]{srgw}\\
\hline
	$1024=2^{10}$ & $93, 930$ & VLS, $H_2(2,5)$ & $2^{10}:(L_5(2)\times \Sym(3))$ & \cite[Table~11.8]{srgw}\\
\hline
	$4098=2^{12}$ & $315,3780$ & VLS, $H_4(2,3)$ & $2^{12}:(GL_2(4)\times\GL_3(4)):2$ & Prop.~\ref{p:onedim}\\
\hline
\end{tabular}
\end{table}

\begin{proof} Let $G\leq\AGL_1(q)$ be a one-dimensional affine rank $3$ group. If $G$ is of odd order, then $q\equiv3\pmod4$ and the irreflexive orbitals of $G$ are opposite tournaments; they are called the {\em Paley tournaments}. The automorphism group of a Paley tournament is an index 2 subgroup in the corresponding one-dimensional affine group $\AGL_1(q)$ (see, e.g., \cite[Subsection~9.7]{Jones}). Therefore, $G^{(2)}\leq\AGL_1(q)$ and there are no exceptions in this case.

Thus, we may assume that $G$ is of even order and, up to a complement, the corresponding rank $3$ graph $\Gamma$  is either a
Van Lint--Schrijver, Paley or Peisert graph, see Theorem~\ref{Muz}. The groups $\Aut(\Gamma)\cap\AGL_1(q)=G^{(2)}\cap\AGL_1(q)$ are described in Lemmas~\ref{VLS} and~\ref{Peisert}. According to Theorem~\ref{a1inv}, the parameters of all possible groups $G$, where $G^{(2)}\neq G^{(2)}\cap\AGL_1(q)$, are listed in Table~\ref{a1exctab}. We consider them case by case.

Let $q=9=3^2$. The both subdegrees of $G$ are equal to $4$, see Table~\ref{a1exctab}. According to \cite[Table 11.8]{srgw}, there is a unique rank $3$ graph $\Gamma$ on $9$ vertices  and its full automorphism group is $3^2:G_0$, where the zero stabilizer $G_0\simeq D_8$ is the dihedral group of order $8$. By Lemmas~\ref{VLS} and~\ref{Peisert}, there are two maximal subgroups of $\GamL_1(3^2)=\langle\widehat\omega,\widehat\phi\rangle$ with the orbits of lengths $4$ on $F^*$. Namely, $H_1=\langle\widehat\omega^2,\widehat\phi\rangle$ and $H_2=\langle\widehat\omega^4,\widehat\phi\widehat\omega\rangle$, the corresponding graphs are the Paley and the Peisert graphs, which are isomorphic in view of Lemma~\ref{PP}. Since $H_1\simeq D_8$, in the first case, $G_1=G_1^{(2)}$, where $G_1=F:H_1$. In the second case, $H_2$ is cyclic of order $4$, and for $G_2=F:H_2$, we have $G_2^{(2)}\not\leq F:\langle\widehat\omega,\widehat\phi\rangle$. However, as it is readily seen, $H_2$ lies in some conjugate of $\langle\widehat\omega,\widehat\phi\rangle$ in $\GL_2(3)$, so we do not have an exception even in this case.

Let $q=16=2^4$. The subdegrees are 5 and 10. By Lemma~\ref{VLS}, the corresponding graph $\Gamma$ is the Van Lint--Schrijver (VLS, for brevity) graph. Then $e=3$ and $|G_0|=|\GamL(q)|/e=20$.  According to \cite[Table 11.8]{srgw}, $G^{(2)}=2^4:\Sym(5)$. Hence $G^{(2)}\not\leq\AGL_1(2^4)$, since $|G_0|<|\Sym(5)|$. Thus, we add $G^{(2)}$ to Table~\ref{not1dim}. It is noted in \cite[Table~11.8]{srgw} that $\Gamma$ is isomorphic to the affine polar graph $\VO^-_4(2)$, so the result also follows from Theorem~\ref{class}(iv)(c), because $\GamO^-_4(2)\simeq\Sym(5)$.

Let $q=25=5^2$. According to Table~\ref{a1exctab}, we are interested in the case when subdegrees are $8$ and $16$. Then, by Lemma~\ref{VLS}, $|G_0|=|\GamL(q)|/3=16$. According to \cite[Table 11.8]{srgw}, $G<G^{(2)}=(\Sym(5)\times\Sym(5)):2$, because $16=|G_0|<|(\Sym(4)\times\Sym(4)):2|=1152$. Therefore, we add the group $G^{(2)}$ to Table~\ref{not1dim}. The corresponding VLS graph is isomorphic to the Hamming graph $H(2,5)=5\times 5$, so the automorphism group preserves the product decomposition $\Omega=\Delta^2$, where $|\Delta|=5$, cf. Theorem~\ref{class}(ii).

Let $q=49=7^2$. In view of Table~\ref{a1exctab}, we need to deal with the case when the subdegrees are equal to $24$. Up to  isomorphism, there are two rank 3 graphs with these parameters: the Paley and Peisert graphs. By \cite[\S~10.18]{srgw}, if the corresponding graph $\Gamma$ is Paley, then $G^{(2)}=7^2:(3\times D_{16})\leq\AGL_1(7^2)$, while if $\Gamma$ is Peisert, then $G^{(2)}=7^2:(3\times\SL_2(3))$ and we add the group $G^{(2)}$ to Table~\ref{not1dim}. The later graph is also an orbital graph of a group from `extraspecial' class (B), see Table~\ref{exceptB}.

Let $q=64=2^6$. The subdegrees are $21$ and $42$. According to \cite[Table 11.8]{srgw}, there is a unique rank 3 graph $\Gamma$ with these parameters and $\Aut(\Gamma)=G^{(2)}=2^6:(\Sym(3)\times L_3(2))$. This group is nonsolvable, so $G^{(2)}=\Aut(\Gamma)$ must be added to Table~\ref{not1dim}. It is worth mentioning that the VLS graph $\Gamma$ is isomorphic to the bilinear forms graph $H_2(2,3)$, so the result also follows from Theorem~\ref{class}(iv)(c).

Let $q=81=3^4$. There are two sets of possible subdegrees here. First, let they be $16$ and $64$, so $\Gamma$ is the VLS graph.  By \cite[Table~11.8]{srgw}, the full automorphism group of $\Gamma$ is $(\Sym(9)\times\Sym(9)):2$, because $\Gamma$ is the Hamming graph $H(2,9)$ preserving the product decomposition $9\times 9$, cf. Theorem~\ref{class}(ii). Thus, $G^{(2)}$ appears in Table~\ref{not1dim}.

The second possibility arises when both subdegrees are of size $40$. Then, see \cite[\S~10.30]{srgw}, there are two rank $3$ graphs $\Gamma$ with these parameters: Paley and Peisert. The full automorphism group of the Paley graph lies in $\AGL_1(3^4)$. If $\Gamma$ is the Peisert graph, then $\Aut(\Gamma)\simeq3^4:(\SL_2(5):2^2)\not\leq\AGL_1(3^4)$, because $\Aut(\Gamma)$ is nonsolvable; and we add the group $G^{(2)}$ to Table~\ref{not1dim}. In fact, the Peisert graph $\Gamma$ can be considered as an orbital graph of a group from `exceptional' class (C) in this case, see Table~\ref{exceptC}.

Let $q=256=2^8$. The subdegrees are 51 and 204. According to \cite[Table 11.8]{srgw}, there is a unique graph $\Gamma$ with these parameters and $\Aut(\Gamma)\simeq2^8:(3\times \SL_2(2^4)):4$, so we add the group $G^{(2)}$ to Table~\ref{not1dim}. The graph $\Gamma$ is also the affine polar graph $\VO_4^-(4)$, so the result follows from Theorem~\ref{class}(iv)(c), because $\GamO^{-}_{4}(4)\simeq(3\times\SL_2(2^4)):4$.

Let $q=289=17^2$. The orbits are of sizes $96$ and $192$. There are, up to isomorphism, two rank $3$ graphs $\Gamma$ with these parameters, see \cite[Table~11.8]{srgw}. The first of them is the VLS graph, because $\Aut(\Gamma)$ is exactly the group described in Lemma~\ref{VLS}, so we do not have an exception in this case. The second graph is an orbital graph corresponding to a group from `extraspecial' class (B), see Table~\ref{exceptB}.

Let $q=529=23^2$. The subdegrees are equal to 264. There are exactly three rank $3$ graphs with these parameters \cite[\S~10.70]{srgw}. Two of them, the Paley and Peisert graphs, are one-dimensional affine, while the third one is the sporadic Peisert which can also be considered as an orbital graph from class (B). Lemma~\ref{PP} implies that there are no exceptions in this case.

Let $q=729=3^6$. The orbits are of sizes $104$ and $624$. By~\cite[Table~11.8]{srgw}, there is a unique rank $3$ graph $\Gamma$ with these parameters and $\Aut(\Gamma)\simeq3^6:(L_3(3)\times\GL_2(3))$. This adds a new row to our table of exceptions. The graph $\Gamma$ is also the bilinear forms graph $H_3(2,3)$, so the result also follows from Theorem~\ref{class}(iv)(b).

Let $q=1024=2^{10}$. The subdegrees are $93$ and $930$. By \cite[Table~11.8]{srgw}, there is a unique rank $3$ graph $\Gamma$ with these parameters and $\Aut(\Gamma)\simeq2^{10}:(L_5(2)\times\Sym(3))$. Therefore, the group $G^{(2)}=\Aut(\Gamma)$ must be added to Table~\ref{not1dim}. The graph $\Gamma$ can be also considered as the bilinear forms graph $H_2(2,5)$, see Theorem~\ref{class}(iv)(b).

Let $q=2209=47^2$. Then two orbits are of the same size $1104$. Therefore, we have the Paley and the Peisert graphs here and do not have exceptions in view of Lemma~\ref{PP}.

Let $q=2401=7^4$.  The subdegrees are $480$ and $1920$. In view of \cite[\S~10.89B]{srgw}, there are two rank 3 graphs $\Gamma$ in this case. If $\Gamma$ is the VLS graph, then $\Aut(\Gamma)\simeq7^4:(480:4)\leq\AGL_1(7^4)$, so no exceptions appear. The second graph is an orbital graph of a rank $3$ group from class (B), see Table~\ref{exceptB}.

Let $q=4096=2^{12}$. By Lemma~\ref{VLS}, the parameter $e$ must be equal to $13$, because the subdegrees are equal to $(q-1)/e=315$ and $(e-1)(q-1)/e=3780$, so $G\simeq2^{12}:(315:12)$. It follows from the proof of \cite[Lemma~3.5]{SkrARS} that if a rank 3 group has same subdegrees, then it induces the bilinear forms graph $H_4(2,3)$. By Theorem~\ref{class}(iv)(b), $\Aut(H_4(2,3)) = GF(4)^{6}:((\GL_2(4) \circ \GL_3(4)) : \Aut(GF(4))).$
Calculations in GAP \cite{gap} show that $|G^{(2)}|=|\Aut(H_4(2,3))|$. Thus, $G^{(2)}=\Aut(H_4(2,3))$ and we add it to the table.

Let $q=7921=89^2$. The subdegrees are $2640$ and $5280$. By Lemma~\ref{VLS}, $e=3$, and $G\simeq89^2:(2640:2)$, where the zero stabilizer $G_0\simeq2640:2$ is of order $5280$. By Table~\ref{a1exctab}, a rank 3 group $H\not\leq\AGL_1(89^2)$ with the same parameters belongs to class (C). Then, according to \cite{liebeckAffine} (see also \cite[Table~8.2]{bray}), the zero stabilizer $H_0$ of $H$ is isomorphic to the central product $F^*\circ\SL_2(5)$ and of order $5280$. Since $|H_0|=|G_0|$ and $H_0$ is nonsolvable, the corresponding graphs cannot be isomorphic, so there is no exception in this case.

Let $q=65536=2^{16}$. Then according to Table~\ref{a1exctab}, the size of the smaller of two orbits must be $3855$. However, $3855=(q-1)/e$, so $e=17$ which contradicts the fact that $\ord_{17}(2)=8\neq16$ (cf. Lemma~\ref{VLS}). Thus, this case is simply impossible.
\end{proof}

\begin{proposition}\label{p:exceptB}
Let  $G$ be one of the rank $3$ groups from the class {\rm(B)}, then $G^{(2)}$ is listed in Table~{\rm\ref{exceptB}}.
\end{proposition}

\begin{table}[ht]
\caption{Closures of groups from class (B)}\label{exceptB}
\begin{tabular}{l|l|l|l|l}
\hline
$q=p^d$ & subdegs & $R$  & $G^{(2)}$ & refs\\
\hline
	$49=7^2$ & $24,24$ & $R_1^1$, $R_2^1$  & $7^2:(3\times SL_2(3))$  & Prop.~\ref{p:onedim} and \cite[\S~10.18]{srgw}\\
\hline
    $64=2^6$ & 27,36 & $3^{1+2}$ &  $2^6:\GamO_6^-(2)$ &  \cite[\S~10.25]{srgw}\\
\hline
	$81=3^4$ & $16,64$ & $R_2^2$   & $(\Sym(9)\times\Sym(9)):2$  & Prop.~\ref{p:onedim} and \cite[Table~11.8]{srgw}\\
\hline
	$81=3^4$ & $32,48$ & $R_1^1$, $R_2^1$   & $3^4:\GamO^+_4(3)$  & Prop.~\ref{p:exceptB} and \cite[Table~11.8]{srgw}\\
\hline
	$81=3^4$ & $32,48$ & $R_1^2$, $R_2^2$   & $3^4:\GamO^+_4(3)$  & Prop.~\ref{p:exceptB} and \cite[Table~11.8]{srgw}\\
\hline
	$169=13^2$ & $72,96$ & $R_1^1$, $R_2^1$ & $13^2:(3\times (SL_2(3):4))$  & \cite[Table~11.8]{srgw}\\
\hline
	$289=17^2$ & $96,192$ & $R_1^1$, $R_2^1$ & $17^2:(8.\Sym(4):2)$  & Prop.~\ref{p:onedim} and \cite[Table~11.8]{srgw}\\
\hline
	$361=19^2$ & $144,216$ & $R_1^1$, $R_2^1$ & $19^2:(9\times\GL_2(3))$  &  \cite[Table~11.8]{srgw}\\
\hline
	$529=23^2$ & $264,264$ & $R_1^1$, $R_2^1$ & $23^2:(11\times\SL_2(3))$  & Prop.~\ref{p:onedim} and \cite[\S~10.70]{srgw}\\
\hline
    $625=5^4$ & $240,384$ & $R_2^2$, $R_3^2$ & $5^4:4.(2^4:\Sym(6))$ & \cite[\S~10.73B]{srgw}\\
\hline
	$729=3^6$ & $104,624$ & $R_1^1$, $R_2^1$ & $3^6:(L_3(3)\times\GL_2(3))$  & Prop.~\ref{p:onedim} and \cite[Table~11.8]{srgw}\\
\hline
	$841=29^2$ & $168,672$ & $R_1^1$, $R_2^1$ & $29^2:(7\times(\SL_2(3):4))$  &  \cite[Table~11.8]{srgw}\\
\hline
	$961=31^2$ & $240,720$ & $R_1^1$, $R_2^1$ & $31^2:(15\times(2.\Sym(4))$  &  \cite[Table~11.8]{srgw}\\
\hline
	$2209=47^2$ & $1104,1104$ & $R_1^1$, $R_2^1$ & $23\times GL_2(3)$ & Prop.~\ref{p:exceptB}\\
\hline
    $2401=7^4$ & $480,1920$ & $R_2^2$ & $7^4:6.(2^{4}:\Sym(5))$ & \cite[\S~10.89B]{srgw}\\
\hline
    $6561=3^8$ & $1440,5120$ & $R_2^3$ & $3^8:(2.2^6:O_6^-(2).2)$ & \cite[\S~10.93]{srgw}\\
\hline
\end{tabular}
\end{table}

\begin{proof} Recall that in this `extraspecial' case, the zero stabilizer $G_0$ of an affine group $G=V:G_0$ normalizes some `extraspecial' subgroup $R$ in $GL(V)$, see Section~\ref{s:known}. Further, we move case by case as in the proof of the previous proposition. The information on the corresponding subgroups $R$ is taken from~\cite[Table~11.5]{srgw}.

Let $q=49=7^2$. Then $R$ is one of the groups $R_1^1$ or $R_2^1$, i.\,e. the dihedral group $D_8$ or quaternion group $Q_8$. The subdegrees are equal to $24$. We have already considered this case in Proposition~\ref{p:onedim} using information from \cite[\S~10.18]{srgw}.

Let  $q=64=2^6$. Then the subdegrees are $27,36$ and $R\simeq3^{1+2}$. In view of \cite[\S~10.25]{srgw}, a unique rank 3 graph $\Gamma$ with these parameters is the affine polar graph $\VO_6^-(2)$. Thus, $\Aut(\Gamma)\simeq2^6:\GamO_6^-(2)$, cf. Theorem~\ref{class}(iv)(c).

Let $q=81=3^4$. There are two sets of possible subdegrees here. First one, when they are $16,64$, was considered in Proposition~\ref{p:onedim}. In the second case, the subdegrees are $32$ and $48$. According to \cite[Table 11.8]{srgw}, there is a unique graph with these parameters and it is isomorphic to the affine polar graph $\VO^+_4(3)$, so the full automorphism group is $3^4:\GamO^+_4(3)\simeq 3^4:((2^{3+2}:3^2):D_8)$, see Theorem~\ref{class}(iv)(c) and \cite[Table 11.8]{srgw}. It is worth noting that a rank 3 group with these parameters appears not only in normalizers of $R_1^1$ and $R_2^1$, but also in normalizers of $R_1^2$ and $R_2^2$, see \cite[Table~11.5 and Theorem~11.4.4]{srgw}, so we add two rows in our table.

Let $q=169=13^2$. Then subdegrees are $72,96$, and by \cite[Table~11.8]{srgw}, there is a unique rank 3 graph $\Gamma$ with these parameters and the automorphism group $\Aut(\Gamma)=V:N_{GL(V)}(R)\simeq13^2:(3\times (SL_2(3):4))$, where $R$ is either $R_1^1$ or $R_2^1$.

In the cases, when $q$ is one of the numbers $361=19^2$, $841=29^2$, and $961=31^2$, the arguments are absolutely the same as in the case $q=169=13^2$, see \cite[Table~11.8]{srgw}.

The cases, when $q$ is one of the numbers $289=17^2$, $529=23^2$, and $729=3^6$, have already appeared in the proof of Proposition~\ref{p:onedim}. When $q\in\{289,529\}$, we have $\Aut(\Gamma)=V:N_{GL(V)}(R)$, where $R$ is $R_1^1$ or $R_2^1$. If $q=3^6$, then the corresponding graph is the bilinear forms graph $H_3(2,3)$ and $\Aut(\Gamma)\simeq3^6:(L_3(3)\times\GL_2(3))$ in view of Theorem~\ref{class}(iv)(b).

Let $q=625=5^4$. The subdegrees are $240,384$, and $R$ is either $R_2^2$ or $R_3^2$.  By \cite[\S~10.73B]{srgw}, there is a unique rank 3 graph $\Gamma$ with these parameters and the automorphism group $\Aut(\Gamma)=V:N_{GL(V)}(R)\simeq5^4:4.(2^4:\Sym(6))$.

Let $q=2209=47^2$.  Then subdegrees are equal to~$1104$. It follows from \cite[Lemma~3.5]{SkrARS} that a group with the same parameters is either from class (B) or is a subgroup in $\AGL_1(q)$. Since $G_0$ includes section isomorphic to $\Sym(4)$, it is not metacyclic, so the corresponding graph is not the Paley or Peisert graph. Thus, $G^{(2)}=N_{GL_2(47)}(R)=23\times GL_2(3)$ in view of~\cite[Theorem~11.4.4]{srgw}.

Let $q=2401=7^4$. The subdegrees are $480,1920$, and $R=R_2^2$.  By \cite[\S~10.89B]{srgw}, there is a unique rank 3 graph $\Gamma$ with these parameters and the automorphism group $\Aut(\Gamma)=V:N_{GL(V)}(R)\simeq7^4:6.(2^4:\Sym(5))$.

Let $q=6561=3^8$. The subdegrees are $1440$, $5120$. In view of \cite[\S~10.93]{srgw}, there is a unique strongly regular graph $\Gamma$ with these parameters and $\Aut(\Gamma)=3^8:(2_-^{1+6}:O_6^-(2).2)$. \end{proof}

\begin{proposition}\label{p:exceptC}
Let  $G$ be one of the rank $3$ groups from the class {\rm(C)}, then $G^{(2)}$ is listed in Table~{\rm\ref{exceptC}}.
\end{proposition}

\begin{table}[ht]
\caption{Closures of groups from (C)}\label{exceptC}
\begin{tabular}{l|l|l|l|l}
\hline
$q=p^d$ & subdegs & $L$ & $G^{(2)}$ & refs\\
\hline
    $81=3^4$ & $40,40$ & $\Alt(5)$ & $3^4:(SL_2(5):2^2)$ & \cite[\S~10.30]{srgw}\\
\hline
    $961=31^2$ & $360,600$ & $\Alt(5)$ & $31^2:(15\times SL_2(5))$ & \cite[\S~10.77]{srgw}\\
\hline
    $1681=41^2$ & $480,1200$ & $\Alt(5)$ & $41^2:(40\circ\SL_2(5))$ & Prop.~\ref{p:exceptC}\\
\hline
    $2401=7^4$ & $960,1440$ & $\Alt(5)$ & $7^4:(48\circ\SL_2(5))$ & \cite[\S~10.89D]{srgw}\\
\hline
    $5041=71^2$ & $840,4200$ & $\Alt(5)$ & $79^2:(35\times\SL_2(5))$ & Prop.~\ref{p:exceptC}\\
\hline
    $6241=79^2$ & $1560,4680$ & $\Alt(5)$ & $79^2:(39\times\SL_2(5))$ & Prop.~\ref{p:exceptC}\\
\hline
	$7921=89^2$ & $2640,5280$ & $\Alt(5)$ & $89^2:(88\circ\SL_2(5))$ & Prop.~\ref{p:exceptC}\\
\hline
	$64=2^6$ & $18,45$ & $\Alt(6)$ & $2^6:3.\Sym(6)$ & \cite[\S~10.24]{srgw}\\
\hline
	$625=5^4$ & $144,480$ & $\Alt(6)$ & $5^4:4.\Sym(6)$ & \cite[\S~10.73A]{srgw}\\
\hline
	$243=3^5$ & $22,220$ & $M_{11}$ & $3^5:(2\times M_{11})$ & \cite[\S~10.55]{srgw}\\
\hline
	$243=3^5$ & $110,132$ & $M_{11}$ & $3^5:(2\times M_{11})$ & \cite[\S~10.55]{srgw}\\
\hline
	$2401=7^4$ & $720,1680$ & $\Alt(7)$ & $7^4:6.\Sym(7)$ & \cite[\S~10.89C]{srgw}\\
\hline
	$256=2^8$ & $45,210$ & $\Alt(7)<\Alt(8)$ & $2^8:(\Alt(8)\times \Sym(3))$ & \cite[\S~10.57]{srgw}\\
\hline
	$256=2^8$ & $120,135$ & $\Alt(9)$ & $2^8:SO_8^+(2)$ & \cite[\S~10.60]{srgw}\\
\hline
	$256=2^8$ & $45,210$ & $\Alt(10)$ & $2^8:\Sym(10)$ & \cite[\S~10.57]{srgw}\\
\hline
	$256=2^8$ & $102,153$ & $L_2(17)$ & $2^8:L_2(17)$ & \cite[\S~10.58]{srgw}\\
\hline
	$729=3^6$ & $224,504$ & $L_3(4), U_3(3)<U_4(3)$ & $3^6:\GamO_6^-(3)$ & \cite[\S~10.76]{srgw}\\
\hline
	$2048=2^{11}$ & $276,1771$ & $M_{24}$ & $2^{11}:M_{24}$ & \cite[\S~10.84]{srgw}\\
\hline
	$2048=2^{11}$ & $759,1288$ & $M_{24}$ & $2^{11}:M_{24}$ & \cite[\S~10.85]{srgw}\\
\hline
	$2401=7^4$ & $240,2160$ & $S_4(3)$ & $7^4:6.S_4(3)$ & \cite[\S~10.89A]{srgw}\\
\hline
	$4096=2^{12}$  & $1575,2520$ & $HJ$ & $2^{12}:(3\times HJ):2$ & \cite[\S~10.92]{srgw}\\
\hline
	$15625=5^{6}$ & $7560,8064$ & $HJ$ & $5^{6}:(2.HJ):4$ & \cite[\S~10.95]{srgw}\\
\hline
	$531441=3^{12}$ & $65520,$  & $G_2(4)<\operatorname{Suz}$ & $3^{12}:2.\operatorname{Suz}.2$ & \cite[\S~6.3.3]{srgw},\\
 & $465920$ & & &  \cite[\S~10.100]{srgw}\\
\hline
\end{tabular}
\end{table}

\begin{proof} For an affine group $G=V:G_0$ from class (C), the socle $L=\Soc(G_0/Z(G_0))$ of the quotient of the zero stabilizer $G_0$ by its center is a finite simple group. If $L\neq\Alt(5)$, then one can find the $2$-closure of $G$ in the corresponding section of \cite[Chapter~10]{srgw} and we add the reference to Table~\ref{exceptC}. Thus, further, we suppose that $L=\Alt(5)$.

If $q$ is the one of the numbers $81=3^4$, $961=31^2$, and $2401=7^4$, then again the $2$-closure of $G$ can be extracted from the corresponding section of \cite[Chapter~10]{srgw}. So we need only to deal with the four remaining cases, when $q$ is one of the numbers $41^2$, $71^2$, $79^2$ and $89^2$. In fact, we will prove that in all these cases $G^{(2)}=V:H_0$, where $H_0=Z(GL_2(q))\circ S$ is a central product of the center of $GL_2(q)$ and a subgroup $S$ of $SL_2(q)$, isomorphic to $SL_2(5)$. From \cite[Theorem 11.4.3]{srgw}, it follows that $G^{(2)}\cap GL(V)=V:H_0$, in particular, we may suppose that $G=V:H_0$. The question that we need to address is whether $G^{(2)}\leq V:GL(V)$ or, in other words, whether $G^{(2)}$ belongs to class~(C)?

If $q=7921=89^2$, then the subdegrees of $G$ coincide with the subdegrees of the one-dimensional affine rank 3 group $A$ of the same degree, see the proof of Proposition~\ref{p:onedim}. Arguing as in this proof, we see that if $G^{(2)}$ is not from class (C), then $G^{(2)}=A$, but it is impossible, since $A$ is solvable. Thus,  $G^{(2)}$ belongs to class (C), as required.

In three remaining cases, we need to check whether there is a rank 3 group that has the same subdegrees as $H=V:H_0$. Suppose first that $G$ and $G^{(2)}$ have the different socles. Since $G$ is not almost simple, \cite[Lemma~2.5]{SkrARS} implies that $G^{(2)}$ and so $G$ must preserve a nontrivial product decomposition. It follows that group $G$ of degree $q$ must have subdegrees $2(\sqrt{q}-1)$ and $(\sqrt{q}-1)^2$. Comparing these parameters with the degrees and subdegrees of the remaining groups $G$, we arrive at contradiction. Thus, the $2$-closure of $G$ must be an affine rank 3 group. We have already considered all the possible coincidences of the parameters for groups from classes (A1) and (B) in Propositions~\ref{p:onedim} and~\ref{p:exceptB}, where (A1) is the subclass of class (A) consisting of the one-dimensional affine rank 3 groups. Thus, we need to compare the parameters of the remaining groups from class (C) with the degrees and subdegrees of the groups from class (A), except subclass (A1). The later parameters are given in~\cite[Table~12]{liebeckAffine}, see also \cite[Table~1]{SkrARS}. Since the degrees in question are the squares of primes, it follows that all the subclasses except (A2) are simply impossible. The subdegrees in the case of subclass (A2) are  $2(\sqrt{q}-1)$ and $(\sqrt{q}-1)^2$, this is also impossible, see the above case of the groups with nontrivial product decomposition. This completes the proof of Proposition~\ref{p:exceptC}, and so the proof of Theorem~\ref{t:2clos}. \end{proof}
\medskip

{\bf Acknowledgments.}  Jin Guo was supported by the Hainan Provincial Natural Science Foundation (Grant no.~122RC543), Jin Guo and Rui Wang were supported by the National Natural Science Foundation of China (No.~12361003), A. V. Vasil'ev was supported by the Sobolev Institute of Mathematics state contract (project FWNF-2022-0002) and by the National Natural Science Foundation of China (No.~12171126).

The authors thank Saveliy Skresanov and Alexander Buturlakin for helpful suggestions, improving the text of the paper.

\end{document}